\documentclass[12pt]{amsart}
\usepackage[T2A]{fontenc}
\usepackage[cp1251]{inputenc}
\usepackage[english]{babel}
\usepackage{amsmath,latexsym,amsthm,amsfonts,tipa,upgreek}
\usepackage{amssymb,amscd,graphpap, stmaryrd}

\newcommand\ZZ{Z}
\newcommand\w{{\omega}}

\newcommand\kk{{\kappa}}

\newcommand\B{{\mathcal B}}
\newcommand\II{{\mathcal I}}

\newtheorem{Th}{Theorem}

\theoremstyle{definition}

\begin{document}

\title{On asymorphisms of groups}
\author{Igor Protasov, Serhii Slobodianiuk}
\date{}
\address{Department of Cybernetics, Kyiv University, Volodymyrska 64, 01033, Kyiv, Ukraine}
\email{i.v.protasov@gmail.com; }
\address{Department of Mechanics and Mathematics, Kyiv University, Volodymyrska 64, 01033, Kyiv, Ukraine}
\email{slobodianiuk@yandex.ru}
\maketitle

\begin{abstract} 
Let $G$, $H$ be groups and $\kk$ be a cardinal. A bijection $f:G\to H$ is caled on asymorphism if, for any $X\in[G]^{<\kk}$, $Y\in[H]^{<\kk}$, there exist $X'\in[G]^{<\kk}$, $Y'\in[H]^{<\kk}$ such that for all $x\in G$ and $y\in H$, we have $f(Xx)\subseteq Y'f(x)$, $f^{-1}(Yy)\subseteq X'f^{-1}(y)$. For a set $S$, $[S]^{<\kk}$ denotes the set $\{S'\subseteq S: |S'|<\kk\}$.

Let $\kk$ and $\gamma$ be cardinals such that $\aleph_0<\kk\le\gamma$.
We prove that any two Abelian groups of cardinality $\gamma$ are $\kk$-asymorphic, but the free group of rank $\gamma$ is not $\kk$-asymorphic to an Abelian group provided that either $\kk<\gamma$ or $\kk=\gamma$ and $\kk$ is a singular cardinal. 
It is known \cite{b6} that if $\gamma = \kk$ and $\kk$ is regular then any two groups of cardinality $\kk$ are $\kk$-asymorphic.
\end{abstract}

\section{Introduction}
Following \cite{b4}, \cite{b7}, we say that a {\em ball structure} is a triple $\B=(X,P,B)$, where $X$, $P$ are non-empty sets, and for all $x\in X$ and $\alpha\in P$, $B(x, \alpha)$ is a subset of $X$ which is called a {\em ball of radius} $\alpha$ around $x$. It is supposed that $x\in B(x, \alpha)$ for all $x\in X$, $\alpha\in P$.
The set $X$ is called the {\em support} of $\B$, $P$  is called the {\em set of radii}.

Given any $x\in X$, $A\subseteq X$, $\alpha\in P$, we set
$$B^*(x,\alpha)=\{y\in X:x\in B(y,\alpha)\},\ B(A,\alpha)=\bigcup_{a\in A}B(a,\alpha),\ B^*(A,\alpha)=\bigcup_{a\in A}B^*(a,\alpha).$$

A ball structure $\mathcal{B}=(X,P,B)$ is called a {\it ballean} if

\begin{itemize}
\item{} for any $\alpha,\beta\in P$, there exist $\alpha',\beta'$ such that, for every $x\in X$,
$$B(x,\alpha)\subseteq B^*(x,\alpha'),\ B^*(x,\beta)\subseteq B(x,\beta');$$

\item{} for any $\alpha,\beta\in P$, there exists $\gamma\in P$ such that, for every $x\in X$,
$$B(B(x,\alpha),\beta)\subseteq B(x,\gamma);$$

\item{} for any $x,y\in X$, there exists $\alpha\in P$ such that $y\in B(x, \alpha)$.
\end{itemize}

We note that a ballean can be considered as an asymptotic counterpart of a uniform space, and could be defined \cite{b8} in terms of the entourages of the diagonal $\Delta_X$ in $X\times X$. In this case a ballean is called a {\em coarse structure}. For categorical look at the balleans and coarse structures as "two faces of the same coin" see \cite{b2}.

Let $\mathcal{B}=(X,P,B)$, $\mathcal{B'}=(X',P',B')$ be balleans. 
A mapping $f:X\to X'$ is called a $\prec$-\emph{mapping} if, for every $\alpha\in P$, there exists $\alpha'\in P'$ such that, for every $x\in X$, $f(B(x,\alpha))\subseteq B'(f(x),\alpha')$. 

A bijection $f:X\rightarrow X'$ is called an {\it asymorphism} between $\B$ and $\B'$ if $f$ and $f^{-1}$ are $\prec$-mappings. In this case $\B$ and $\B'$ are called {\it asymorphic}.

Let $\B = (X,P,B)$ be a ballean. Each subset $Y$ of $X$ defines a {\it subballean} $\B_Y = (Y,P,B_Y)$, where $B_Y(y,\alpha) = Y \cap B(y, \alpha)$. A subset $Y$ of $X$ is called {\it large} if $X = B(Y, \alpha)$, for $\alpha \in P$. Two balleans $\B$ and $\B'$ with supports $X$ and $X'$ are called {\it coarsely equivalent} if there exist large subsets $Y\subseteq X$ and $Y' \subseteq X'$ such that the subballeans $\B_Y$ and $\B'_{Y'}$ are asymorphic.
In the proof of Theorem~\ref{t2}, we use the following equivalent definition: $\B$ and $\B'$ are coarsely equivalent if there is a $\prec$-mapping $f:X_1\to X_2$ such that $f(X_1)$ is large and, for every $\alpha'\in P'$ there exists $\alpha \in P$ such that $f^{-1}(B'(f(x), \alpha'))\subseteq \B(x, \alpha)$ for each $x\in X$.

We recall \cite{b5} that an ideal $\II$ in the Boolean algebra of all subsets of a group $G$ is a {\em group ideal} if $AB^{-1}\in\II$ for any $A,B\in\II$. We suppose that $\II$ contains all finite subsets of $G$ and denote by $(G,\II)$ the ballean $\B=(G,\II,B)$, where $P(g, A) = (A\cup\{e\})g$, $e$ is the identity of $G$.

For an infinite group $G$ and an infinite cardinal $\kk$, $\kk\le|G|$, we denote by $[G]^{<\kk}$ the group ideal $\{A\subset G: |A| < \kk\}$.

We say that two groups $G$ and $H$ are {\em $\kk$-asymorphic} ({\em $\kk$-coarsely equivalent}) if the balleans $(G, [G]^{<\kk})$ and $(H, [H]^{<\kk})$ are asymorphic (coarsely equivalent).

In the case $\kk=\aleph_0$, we say that $G$ and $H$ are {\em finitary asymorphic} and {\em finitary coarsely equivalent} respectively. We note that finitely generated groups are finitary coarsely equivalent if and only if $G$ and $H$ are quasi-isometric \cite[Chapter 4]{b2.5}.

A classification of countable locally finite groups (each finite subset generates finite subgroup) up to finitary asymorphisms is obtained in \cite{b3}: there are continuum distinct types (see also \cite[p. 103]{b4}) and each such group is finitary asymorphic to some direct product of finite cyclic groups. For coarse classifications of countable Abelian groups see \cite{b1}. Any two countable torsion Abelian groups are finitary coarsely equivalent.

This note is motivated by the following result \cite[Theorem 3]{b6}. Let $G$, $H$ be groups of cardinality $\gamma$, $\gamma>\aleph_0$ and let $\kk$ be a cardinal such that $\aleph_0<\kk\le\gamma$. If $\kk = \gamma$ and $\kk$ is regular then $G$ and $H$ are $\kk$-asymorphic.

What happens if $\kk < \gamma$ or $\kk = \gamma$ but $\kk$ is singular?

\section{Theorems}
\begin{Th}\label{t1} 
Let $G$ be an Abelian group of cardinality $\gamma$, $\gamma>\aleph_0$ and let $\kk$ be a cardinal such that $\aleph_0<\kk\le\gamma$. Then $G$ is $\kk$-asymorphic to the free Abelian group $A_\gamma$ of rank $\gamma$.
\end{Th}
\begin{Th}\label{t2} Let $G$ be an Abelian group of cardinality $\gamma$, $\gamma>\aleph_0$. Then $G$ is not $\kk$-coarsely equivalent to the free group $F_\gamma$ of rank $\gamma$ provided that either $\aleph_0 < \kk < \gamma$ or $\kk= \gamma$ and $\gamma$ is a singular cardinal. In particular, $G$ and $F_\gamma$ are not $\kk$-asymorphic.
\end{Th}

\section{Proofs}
{\em Proof of Theorem~\ref{t1}.}
We choose a system of subgroups $\{ G_\alpha : \alpha < \gamma \}$ of $G$ such that
\begin{itemize}
\item[(1)] $G_0=\{e\}$, $G=\bigcup_{\alpha<\kk}G_\alpha$;
\item[(2)] $G\alpha\subset G_\beta$ for all $\alpha < \beta < \gamma$;
\item[(3)] $G_\beta=\bigcup_{\alpha<\beta} G_\alpha$ for every limit ordinal $\beta < \gamma$;
\item[(4)] $|G_{\alpha+1}:G_\alpha|=\aleph_0$ for every $\alpha < \gamma$.
\end{itemize}
For each $\alpha<\gamma$, we fix some system $X_\alpha$ of representatives of cosets of $G_{\alpha+1}\setminus G_\alpha$ by $G_\alpha$ so $G_{\alpha+1}\setminus G_\alpha = X_\alpha G_\alpha$.

We take an arbitrary element $g\in G\setminus\{e\}$ and choose the smallest subgroup $G_\gamma$ such that $g\in G_\alpha$. By $(3)$, $\alpha = \alpha_0+1$ for some $\alpha_0 < \gamma$. Then $g = g_0x_0$,
$g_0\in G_{\alpha_0}$, $x_0\in X_{\alpha_0}$. If $g_0\neq e$, we repeat the argument for $g_0$: choose $\alpha_1$ such that $g_0\in G_{\alpha_1+1}\setminus G_{\alpha_1}$ and write $g_0 = g_1x_1$, where $g_1\in G_{\alpha_1}$, $x_1\in X_{\alpha_1}$ and so on. Since the set of ordinals less than $\kk$ is well ordered, after finite number of steps, we get
$$g = x_{s(g)}\dots x_1x_0,\ \ x_i\in X_{\alpha_i},\ i\in\{0,\dots,s(g)\},\ \ \alpha_0>\alpha_1>\dots>\alpha_{s(g)},$$

We observe that such a representation is unique and denote
$$supt(g) = \{\alpha_{s(g)},\dots, \alpha_1, \alpha_0\},\ \ supt(e) = \varnothing$$
We identify $A_\gamma$ with the direct product $\times_{\alpha<\gamma} \ZZ_\alpha$ of infinite cyclic groups and, for each $\alpha < \gamma$, fix some bijection $f_\alpha : X_\alpha\to \ZZ_\alpha\setminus\{e_\alpha\}$, $e_\alpha$ is the identity of $\ZZ_\alpha$. We define a bijection $f: G\to A_\gamma$ putting $f(e)=(e_\alpha)_{\alpha<\gamma}$ and, for $g\in G\setminus\{e\}$,
$$f(g)=f(x_{s(g)}\dots x_1 x_0)=(f_{\alpha_{s(g)}}(x_{s(g)}),\dots, f_{\alpha_1}(x_1), f_{\alpha_0}(x_0)).$$
We show that $f$ is an asymorphism between the balleans $(G, [G]^{<\kk})$ and $(A_\gamma, [A_\gamma]^{<\kk})$.

To show that $f^{-1}$ is a $\prec$-mapping, we take $a\in A_\gamma$,  $K\in [A_\gamma]^{<\kk}$ and choose a set $I\in[\gamma]^{<\kk}$ such that $K\subseteq \times_{\alpha\in I}\ZZ_\alpha$.
We denote $X = \times_{\alpha\in I}\ZZ_\alpha$, $b=pr_{I}a$, $c=pr_{\gamma\setminus I}a$.
Then we have $$f^{-1}(Ka)\setminus f^{-1}(Xc)=f^{-1}(X)f^{-1}(c)=f^{-1}(X)(f^{-1}(b))^{-1}f^{-1}(b)f^{-1}(c)=$$
$$ = f^{-1}(X)(f^{-1}(b))^{-1}f^{-1}(a)\subseteq f^{-1}(X)(f^{-1}(X))^{-1}f^{-1}(a)$$
and it is suffices to note that $f^{-1}(X)(f^{-1}(X))^{-1}\in [G]^{<\kk}$.
 
The verification that $f$ is a $\prec$-mapping is more delicate.
We take an arbitrary $F\in [G]^{<\kk}$ and denote by $Y$ the smallest subgroup of $G$ containing $F$ and such that if $g\in Y$ and $\alpha\in supt(g)$ then $X_\alpha\subseteq Y$. We show that $Y\in[G]^{<\kk}$. Indeed, $Y$ can be obtained in the following way. For a subset $A$ of $G$, we denote by $<A>$ the subgroup generated by $A$ and $h(A)=A\cup\bigcup\{X_\alpha:\alpha\in supt(g), g\in A\}$.
We put $S_0=<F>$, $Y_0=h(S_0)$ and inductively $S_{n+1}=<Y_n>$, $Y_{n+1}=h(S_n)$. 
Then $Y=\bigcup_{n\in\w}Y_n$.

To conclude the proof, we take an arbitrary $g\in G$, put $I=\bigcup_{g\in Y}supt(g)$, and write $g=g_0g_1$ where $supt(g_0)\subseteq I$, $supt(g_1)\subseteq \gamma\setminus I$. Then we have
$$f(Kg)\subseteq f(Yg_0g_1)\subseteq f(Yg_1)\subseteq f(Y)f(g_1)\subseteq f(Y)f(g).$$

{\em Proof of Theorem~\ref{t2}.}
We are going to get a contradiction assuming only that there is a $\prec$-mapping $f: G\to F_\gamma$ such that $F_\gamma = Kf(G)$ for some $K\in [F_\gamma]^{<\kk}$.
In view of Theorem~\ref{t1}, we may suppose that $G$ is a group of exponent $2$.
Since $f$ is a $\prec$-mapping, for every $a\in G$, there exists $K_a\in[F_\kk]^{<\kk}$ such that, for each $x\in G$, we have $$f(a+x)\in K_a f(x).$$
We note that the family $\{K_a: a\in G\}$ can be chosen so that, for some cardinal $\delta$, $\kk\le\delta<\gamma$, 
\begin{itemize}
\item[(5)] $|K_a|<\delta$ for each $a\in G$.
\end{itemize}
If $\kk < \gamma$ then $(5)$ is evident with $\gamma = \kk$. We consider the case $\gamma = \kk$ and $\kk$ is singular. If $(5)$ could not be satisfied then, for every $\delta <\gamma$, there exists $a_\delta\in G$ such that
$$|\{f(a_\delta+x)f^{-1}(x) : x\in G\}|>\delta.$$

We denote $Y_\delta = \{f(a_\delta+x)f^{-1}(x) : x\in G\}$ and use singularity of $\gamma$ to choose a subset $\Delta$ of $\gamma$ such that $|\Delta|<\gamma$ and the set $\{|Y_\delta|: \delta\in\Delta\}$
is cofinal in $\gamma$. We put $A=\{a_\delta:\delta\in\Delta\}$. Since $f$ is a $\prec$-mapping, there is a subset $Y$ of $F_\gamma$ such that $|Y| = \kk$ and $\{f(A+x)f^{-1}(x):x\in G\}\subseteq Y$.
Then we have got a contradiction with $Y_\delta \subseteq Y$ for each $\delta\in\Delta$.

We consider $F_\gamma$ as the group of all reduced group words over the alphabet $\kk$. For 
$g\in F_\gamma$, we denote by $alph f(g)$ the set of all letters $\alpha <\gamma$ such that $\alpha$ or $\alpha^{-1}$ occures in $g$, and for subset $S$ of $F_\gamma$, we put $alph(S)=\bigcup_{g\in S}alph(g)$.

Now we show how to find $a, b\in G$ such that 
\begin{itemize}
\item[(6)] $alph f(a) \cap alph f(b) = \varnothing$;
\item[(7)] $alph f(a) \setminus alph K_b \neq \varnothing$
\item[(8)] $alph f(b) \cap alph K_a \neq \varnothing$
\end{itemize}

Since $F_\gamma=Kf(G)$ for some $K\in[F_\gamma]^{<\kk}$, we can choose a subset $A\subseteq G$ such that $|A|=\delta$ and $|alph f(A)|=\delta$. We take
$$\beta\in \gamma\setminus(alph f(A)\cup \bigcup_{a\in A}alph K_a)$$
and find $b\in G$ such that $f(b) = \beta$. Since $|K_\beta|<\delta$ and $|alph f(A)| = \delta$, there is $a\in A$ such that $alph f(a)\setminus alph K_b \neq \varnothing$.

We put $f(a+b) = z$ and note that
\begin{itemize}
\item[(9)] $f(b) = f(a+(a+b))\in K_az$;
\item[(10)] $f(a) = f(b+(a+b))\in K_bz$.
\end{itemize}
We take $u\in alph f(a)\setminus K_b$, $v\in alph f(a)\setminus alph K_b$.
Then $u\in alph(z)$, $v\in alph(z)$. If $u$ occures in $z$ before $v$ then, by $(6)$, $(10)$ is false.
If $v$ occurs in $z$ before $u$ then, by $(6)$, $(9)$ is false. These contradictions conclude the proof.


\begin{thebibliography}{99}

\bibitem{b1} T.~Banakh, J,~Higes, M.~Zarichnyi, {\it The coarse classification of countable abelian groups}, Trans. Amer. Math. Soc. {\bf 362} (2010) 4755-4780.

\bibitem{b2} D.~Dikranjan, N.~Zava, {\it Some categorical aspects of coarse spacaes and balleans}, preprint.

\bibitem{b2.5} P. de la Harpe, {\it Topics in Geometric Group Theory}, University Chicago Press, 2000.

\bibitem{b3} I.V.~Protasov, {\it Morphisms of ball structures of groups and graphs}, Ukr. Mat. Zh. {\bf 53} (2002) 847-855.

\bibitem{b4} I.~Protasov, T.~Banakh, {\it Ball structures and colorings of groups and graphs}, Math. Stud. Monogr. Ser., Vol. 11, VNTL, Lviv, 2003.

\bibitem{b5} I.V.~Protasov, O.I.~Protasova, {\it Sketch of group balleans}, Mat. Stud. {\bf 22} (2004) 10-20.

\bibitem{b6} I.V.~Protasov, A.~Tsvietkova, {\it Decomposition of cellular balleans}, Topology Proc. {\bf 36} (2010) 77-83.

\bibitem{b7} I.~Protasov, M.~Zarichnyi, {\it General Asymptology}, Math. Stud. Monogr. Ser., Vol. 12, VNTL, Lviv, 2007.

\bibitem{b8} Roe~J., {\it Lectures on coarse geometry}, Amer. Math. Soc., Providence, R.I, 2003.
\end{thebibliography}
\end{document}